\input amstex
\input xy
\input epsf
\xyoption{all}
\documentstyle{amsppt}
\document
\magnification=1200
\NoBlackBoxes
\nologo
\hoffset1.5cm
\voffset2cm
\pageheight {16cm}




\centerline{\bf FOUNDATIONS AS SUPERSTRUCTURE\footnotemark1} 
\footnotetext{Talk at the international interdisciplinary conference
``Philosophy, Mathematics, Linguistics: Aspects of Interaction'' , Euler Mathematical
Institute, May 22--25, 2012.}
\smallskip

\centerline {\it (Reflections of a practicing mathematician)}

\bigskip

\centerline{\bf Yuri I. Manin}

\medskip

\centerline{\it Max--Planck--Institut f\"ur Mathematik, Bonn, Germany,}

\bigskip

{\it ABSTRACT.} This talk presents foundations of mathematics as a historically variable
set of principles appealing to various modes of human intuition and devoid of any
prescriptive/prohibitive power. At each turn of history, foundations crystallize the accepted norms
of interpersonal and intergenerational transfer and justification of
mathematical knowledge.

\bigskip

\centerline{\bf Introduction}

\medskip

{\bf Foundations vs Metamathematics.} In this talk, I will interpret  the idea of Foundations in the wide sense. For me,
Foundations at each turn of history embody currently recognized, but historically variable, principles of organization of mathematical knowledge and of the
interpersonal/transgenerational transferral of this knowledge. 
When these principles are studied 
using the tools of mathematics itself,  we get a new chapter of mathematics, {\it metamathematics}.
\smallskip

Modern philosophy of mathematics is often preoccupied  with {\it informal interpretations of theorems,
proved in metamathematics of the XX--th century,} of which the most
influential was probably G\"odel's 
incompleteness theorem that aroused considerable existential anxiety.

\smallskip
In metamathematics, G\"odel's theorem is a discovery
that a certain class of finitely generated structures (statements in a formal language)
contains substructures that are 
not finitely generated (those statements that are true in a standard interpretation).

\smallskip

It is no big deal for an algebraist, but certainly interesting thanks to a new context. 

\smallskip

Existential anxiety can be alleviated if one strips 
``Foundations'' from their  rigid {\it prescriptive/prohibitive}, or {\it  normative} functions and  considers
various foundational matters  simply from the viewpoint of their mathematical 
content and on the background of whatever historical period.

\smallskip

Then, say, the structures/categories
controversy is seen in a much more realistic light:
contemporary studies fuse (Bourbaki type) structures and categories freely,
naturally and unavoidably. 

\smallskip
For example, in the definition of {\it abelian categories}
one starts with {\it structurizing sets of morphisms:} they become abelian groups.
In the definition of 2--categories, the sets of morphisms are even {\it categorified:} they become
objects of categories, whose morphisms become then the {\it morphisms of the second level} of initial category.
Since in this way one often obtains vast mental images of complex combinatorial structure,
one applies to them principles 
of {\it homotopy topology} (structural study of topological structures up to homotopy equivalence)  in order to squeeze it down to size etc.
\smallskip

I want to add two more remarks to this personal credo.

\smallskip

First, the recognition of quite restrictive and historically changing
normative function of Foundations makes this word somewhat too
expressive for its content. In a figure of speech such as
``Crisis of Foundations'' it suggests a looming crash of the whole 
building (cf. similar concerns expressed by R.~Hersh, [He]).

\smallskip

But, second, the first ``Crisis of  Foundations'' occurred in a very interesting
historical moment, when the images of formal mathematical reasoning and algorithmic computation
became so precise and detailed that they could be, and were,
described as {\it  new mathematical structures}: formal languages and their interpretations,
partial recursive functions. They could easily fit Bourbaki's universe,
even if Bourbaki himself was too slow and awkward to really appreciate the new development.

\smallskip

At this juncture, contemporary {\it  ``foundations''}  morphed into a {\it superstructure,} 
high level floor of mathematics building itself. This is the reason why I keep using
the suggestive word  ``metamathematics'' for it.

\smallskip
This event  generated a stream of philosophical thought
striving to recover the lost normative function. 
One of the reasons of my private mutiny against it (see e.g.  [Ma1])
was my incapability to find any of the  philosophical arguments more convincing that
even the simplest mathematical reasonings, whatever ``forbidden'' notions they might involve.
\smallskip

In particular, whatever doubts one might have  about the
 scale of Cantorial cardinal and ordinal infinities,
the basic idea of set embodied in Cantor's  famous
``definition'', as a collection of  definite, distinct objects
of our thought, is as alive as ever. Thinking about
a topological space, a category, a homotopy type,
a language or a model, we start with imagining such a collection,
or several ones, and continue by adding new types
``of distinct objects of our thought'', derivable from the previous ones or embodying
fresh insights.
\smallskip

 To summarize: good metamathematics is 
good mathematics rather than shackles on good mathematics.

\medskip

{\bf Plan of the article.} Whatever one's attitude to mathematical Platonism might be,
it is indisputable that human minds constitute an important part of  habitat
of mathematics. In the first
section, I will postulate three basic types of mathematical intuition and argue
that one can recognize them at each scale of study: personal,
interpersonal and historical ones.

\smallskip

The second section is concerned with historical development of the dichotomy
{\it continuous/discrete} and  evolving interrelations between its terms.
\smallskip
Finally, in the third section I briefly recall the
discrete structures of linear languages studied in {\it classical metamathematics,}
and then sketch the growing array of language--like non--discrete structures
that gradually become the subject--matter of {\it contemporary metamathematics.}

\bigskip

\centerline{\bf 1. Modes of mathematical intuition}



\medskip

{\bf 1.1. Three modes.}  I will adopt here the viewpoint according to which {\it at the individual level}
mathematical intuition, both primary and trained one, has three basic sources, that I will
describe as {\it spatial, linguistic,} and {\it operational} ones.

\smallskip

The neurobiological correlates of the {\it spatial/linguistic} dichotomy were elaborated
in the classical studies of {\it lateral asymmetry of brain.} When its mathematical
content is objectivized, one often speaks about the opposition
{\it continuous/discrete}. 
\smallskip

The {\it linguistic/operational}
dichotomy is observed in many experiments studying proto--mathematical abilities
of animals. Animals, when they   solve and communicate solutions of elementary problems 
related to  counting, use not words but {\it actions}:   cf. some expressive
descriptions by Stanislas Dehaine in [De], Chapter 1: ``Talented and gifted animals''.
Operational mode, when it is  externalized and codified, becomes a powerful tool
for social expansion of mathematics. Learning by rote of ``multiplication table'' 
became almost a symbol of democratic education.

\smallskip

The sweeping subdivision of  mathematics into Geometry and  Algebra,
to which at the beginning of modern era was added Analysis (or Calculus)
can be considered as a correlate on the scale of whole (Western) civilization
of  the trichotomy that we postulated above on the scale of an individual (cf. [At]).

\medskip

 It is less widely recognized  that {\it even at the civilization scale}, at various historical
periods, each of the spatial,
linguistic and operational modes of mathematical intuition can dominate and govern the way that
basic mathematical abstractions are perceived and treated. 
\smallskip

I will consider as an example
``natural'' numbers.  Most of us nowadays immediately associate them with their {\it names}:
 decimal notation 1,2,3, \dots, 1984, ... , perhaps completed by  less systemic 
signs such as $10^6$  or $XIX$.

\smallskip

This was decidedly not always so as the following examples stretching over centuries 
and millennia show.

\medskip

 {\bf 1.2.  Euclid and his ``Elements'': spatial and operational {\it vs} linguistic.} For Euclid,
a number was a ``magnitude'', a potential  result of measurement. Measurement of
a geometric figure $A$ by a ``unit'' , another geometric figure $U$,  was conceived as
a ``physical activity in mental space'': moving a segment of line inside another
segment, step by step; paving a square by smaller squares etc. Inequality $A<B$ roughly speaking,
meant that a figure $A$ could be moved to fit inside $B$ (eventually, after cutting $A$  into several pieces and
rearranging them in the interior of $B$).

\smallskip

In this sense, Euclidean geometry might be conceived as ``physics of solid bodies in the
dimensions one, two and three'' (or more precisely, after Einstein,  physics {\it in gravitational vacuum}
of respective dimension).
This pervasive identification of Euclidean space with our physical space probably influenced 
the history of Euclid's  ``fifth postulate''. This history includes repeating attempts {\it to prove it}, that is, to
deduce properties of space ``at infinity''  from observable ones at a finite distance,
and then only  reluctant  accepting  the B\'olyai and Lobachevsky non--Euclidean spaces
as ``non--physical'' ones.

\smallskip
As opposed to addition and subtraction, {\it the multiplication of numbers} naturally required
passage into {\it a higher dimension}: multiplying two lengths, we get a surface. This was
a great obstacle, but, I think, also opened for trained imagination the door
to higher dimensions. At least, when Euclid has to speak about the product of an
arbitrary large finite set of primes (as in his proof involving
$p_1\dots p_n +1$), he is careful to explain his general reasoning
by the case of {\it three} factors, but without doubt, he had some mental images overcoming this
restriction. 

\smallskip
In fact, the strength of spatial and operational imagination required and achieved by modern mathematics
can be glimpsed on a series of examples, starting, say with Morse theory and reaching Perelman's proof
of Poincar\'e conjecture. Moreover, physicists  could
produce such wonders as Feynman's  path integral and Witten's topological invariants,
which mathematicians include in their more rigidly organized world only with considerable
efforts.

\smallskip

At first sight, it might seem strange that the notion of {\it a prime number}, 
theorem about  (potential) {\it infinity of primes}, and theorem about
{\it unique decomposition} could have been stated and proved by Euclid in his geometric
world, when no systematic notation for integers
was accepted as yet, and no computational rules {\it dealing with such a notation} rather than  numbers
themselves were available. 

\smallskip
But trying to rationalize this historical fact, one comes to a somewhat paradoxical 
realization that an efficient notation, such as Hindu--Arabic numerals, actually {\it does not help}, and even 
{\it hinders} the understanding of properties related to divisibility, primality etc. that is, all properties
that refer to numbers themselves rather than their {\it names.}

\medskip

In fact, the whole number theory could come into being only unencumbered by any efficient {\it notation}
for numbers.

 \medskip
 
 {\bf 1.3. ``Algorist and Abacist'':  linguistic {\it vs} operational.} The dissemination of a positional number system in Europe after the
appearance of Leonardo Fibonacci's {\it Liber Abaci} (1202) was, in
essence, the beginning of the expansion of a universal, truly global
language.  Its final victory took quite some time.

\smallskip

The book by Gregorio Reisch, {\it Margarita Philosophica}, was published in Strasbourg in 1504.
One engraving in this book shows a female figure symbolizing Arithmetics.
She contemplates two men, sitting at two different tables,
an {\it abacist} and an {\it algorist}.

\smallskip

The abacist is bent over his {\it abacus}. This primitive calculating device
survived until the days of my youth: every cashier in any shop in
Russia, having accepted a payment, would start calculating change
clicking movable balls of her abacus.

\smallskip

The algorist is computing something, writing Hindu--Arabic numerals  on his desk.
The words  ``algorist'' and modern ``algorithm'' are derived from the name
of the great Al Khwarezmi (born in Khorezm c. 780).

\smallskip

In the context of this subsection, abacus illustrates the operational mode whereas 
computations with numerals do the same for linguistic one (although
in other contexts the operational side of such computations might dominate).
\smallskip
This engraving 
in the reception of contemporary readers was more politicised. It symbolized coming
of a new epoch of democratic learning.

\smallskip

Catholic Church supported the Roman tradition, usage of Roman numerals.
They were fairly useless for practical  commercial bookkeeping,
calender computations such as  dates of Easter and other moveable feasts etc.  
Here abacus  was of great help.

\smallskip

The competing tribe of algorists were able to compute things by writing strange signs on paper or sand,
and their art was associated with dangerous, magical, secret Muslim knowledge.
Al Khwarezmi teaching became their (and our) legacy. 

\smallskip

Arithmetics blesses both practitioners.

 \medskip
 
 {\bf 1.4. John Napier and Alan Turing: operational.}  The nascent programming languages for centuries existed only as  informal 
subdialects of a natural language. They had a very limited (but crucially important) sphere of
applicability, and were addressed to human calculators, not electronic 
or mechanical ones. Even Alan Turing in the 20th century, when speaking of his 
 universal formalization of computability, later called Turing machine, used the word ``computer'' to 
 refer to {\it a person} who mechanically follows a finite  list of instructions lying 
before him/her.
\smallskip

The ninety--page table of natural logarithms that John Napier
published in his book {\it Mirifici Logarithmorum Canonis Descriptio} in 1614
was a paradoxical example of this type of activity that became a
cultural and historical monument on a global scale.
 Napier, who computed the logarithms
manually, digit by digit, combined in one person the role of
creator of new mathematics and that of computer--clerk who
followed his own instructions. His assistant Henry Briggs later performed this
function.

\smallskip

Napier's tables were tables of  (approximate values of) {\it natural logarithms},
with  the  base $e=2,718281828 ...$. However, it seems that he neither
referred to $e$ explicitly, nor even recognized its existence. Roughly speaking,
after having chosen the precision which he wanted to calculate logarithms,
say with error $<10^{-7}$, he dealt with integer powers of the number
$1+10^{-8}$, whose $10^8$ power was close to $e$.

\smallskip

This is one more example of the seemingly paradoxical fact,  that an efficient and unified notation 
for objects of mathematical world can {\it hinder a theoretical understanding} of this world.

\smallskip

All the more amazing was the philosophical insight of Leibniz, who in
his famous exhortation {\it Calculemus!} postulated that not only numerical
manipulations, but any rigorous, logical sequence of thoughts that
derives conclusions from initial axioms can be reduced to
computation.  It was the highest achievement of the great logicians
of the 20th century (Hilbert, Church, G\"odel, Tarski, Turing,
Markov, Kolmogorov,...) to draw a precise map of the boundaries
of the Leibnizian ideal world, in which

\medskip

{\it reasoning is equivalent to computation;

\smallskip

truth can be formalized, but cannot always be verified formally;

\smallskip

 the ``whole truth''  even about the smallest infinite
mathematical universe, natural numbers, exceeds potential of any finitely generated 
language to generate true theorems.} 

\medskip

The central concept of this program, {\it formal languages}, inherited
the basic features of both natural languages (written form fixed by an
alphabet) and the positional number systems of arithmetic.  In particular,
any classical formal language is one--dimensional (linear) and consists of
discrete symbols that explicitly express the basic notions of logic.

\medskip

Euclid found the remedy for the deficiencies of  this linearity by strictly 
restricting role of natural language to the expression
of {\it logic} of his proofs. The {\it content} of his mathematical imagination
was transmitted  by pictures.

 \bigskip
 
 \centerline{\bf  2. Continuous or discrete?}
 
 \smallskip
 
 \centerline{\bf From Euclid to Cantor to homotopy theory }

\medskip

 {\bf 2.1. From continuous to discrete in ``Elements''.} As we have seen, integers (and a restricted amount of other real numbers) 
for Euclid were results of (mental) measurement: {\it discrete came from continuous.}
This was  one--way road: continuous could not be produced from discrete.
The idea that a line ``consists'' of points, so familiar to us today, does not
seem to belong to Euclid's mental world and, in fact, to mental worlds of many
subsequent generations of mathematicians until Georg Cantor. For Euclid,
a point can be (a part of) the boundary of a (segment) of line, but such a segment cannot 
be scattered to a heap of points.

\smallskip

Geometric images are the source and embodiment not only of
numbers, but of logical reasoning as well: in ``Elements'' at least a comparable part of
its logic is encoded in figures rather than in words.

\smallskip

This was made very clear in the London publication of 1847, entitled
\medskip
\centerline{The first six books of}

\centerline{\it THE ELEMENTS OF EUCLID}

\centerline{in which coloured diagrams and symbols} 

\centerline{are used instead of letters for the}

\centerline{greater ease of learners}

\medskip

whose author was Oliver Byrne, {\it ``Surveyor of her Majesty's settlements in the Falklands Islands''.}

\smallskip

Byrnes literally writes algebraic formulas whose main components are triangles, colored sectors of circle,
segments of line etc. connected by more or less conventional algebraic signs.

\medskip

 {\bf 2.2. From discrete to continuous: Cantor, Dedekind, Hausdorff, Bour\-baki ...}
This way is so familiar to my contemporaries that I do not have to spend much time to its
description. The  description of a mathematical structure, such as a group, or a topological space,
 according to Bourbaki starts with one or
several unstructured sets, to which one adds elements of a these sets or derived sets
satisfying restrictions formulated in terms of set theory. 

\smallskip

Thus the twentieth century idea of ``continuous'' is based upon two parallel notions: that of {\it topological space $X$} (a set with the system of ``open'' subsets) and that of a ``continuous map''  $f:\,X\to Y$ between
topological spaces. Further elaboration involving sheaves, topoi etc does not part with this
basic intuition.

\smallskip

However, the set--theoretic point of depart helped enrich the geometric intuition
by images that were totally out of reach  earlier. The  discovery of 
difference between {\it continuous} and {\it measurable} (from Lebesgue integral
to Brownian motion to Feynman integral) was a radical departure from Euclidean universe.
\smallskip
In a finite--dimensional context,
one could now think  about  Cantor sets, Hausdorff dimension and fractals, curves filling a square,
Banach--Tarski theorem. In infinite--dimensional contexts wide new horizons opened, starting with
topologies of Hilbert and Banach linear spaces and widening in an immense universe of topology
and measure theory of non--linear function spaces.

\medskip

{\bf 2.3. From continuous to discrete: homotopy theory.} One of the most important
development of topology was the discovery of main definitions and results of homotopy theory.
Roughly speaking, a homotopy between two topological spaces $X,Y$ is a continuous deformation
producing $Y$ from $X$, and similarly a  homotopy between two continuous maps $f,g:\,X\to Y$
is a continuous deformation producing $g$ from $f$. A {\it homotopy type} is the class of spaces
that are homotopically equivalent pairwise. To see how drastically the homotopy 
can change a space, one can note that a ball, or a cube, of any dimension is {\it
contractible}, that is, can be homotopically deformed
to a point, so that dimension ceases to be invariant.

\smallskip

The basic discrete invariant of the homotopy type of $X$ is the set of
its connected components $\pi_0(X)$. To see, how this invariant gives rise
to one of the basic structures of mathematics, ring of integers $\bold{Z}$,
consider a real plane $P$ with a fixed orientation, a point $x_0$ on it, different from $(0,0)$, and the set
of homotopy classes of loops (closed paths) in $P$, starting and ending at $x_0$ and avoiding
$(0,0)$.  This latter set can be canonically identified with $\bold{Z}$: just count the number
of times the loop goes around $(0,0)$. Each loop going in the direction of
orientation counts as $+1$, where as the ``counter--clockwise'' loops count as $-1$.
\smallskip

On a very primitive level, this identification shows how the ideas of homotopy naturally introduce
negative numbers. In the historically earlier periods when integers were measuring 
geometric figures (or counting real/mental objects) even idea of zero was very difficult and slowly
gained ground in the symbolic framework of positional
notation.
Introduction of negative numbers required appellation to an extra--mathematical reality,
such as {\it debt} in economics.

\smallskip

More generally, Voevodsky in his research project [Vo] introduces 
the following hierarchy of homotopy types graded by their $h$--levels.
Zero level homotopy type consists of one point representing contractible spaces.
If types of level $n$ are already defined, types of level $n+1$ consist of spaces
such that the space of paths between any two points belongs to type of level $n$.

\smallskip

He further interprets types of level 1, represented by one point and empty sets, 
as {\it truth values}, and types of level 2 as sets. All sets in this universe are thus of the form $\pi_0(X)$.

\smallskip

Higher levels are connected with theory of categories, poly--categories etc, and we will return to them
in the next section. At this point, we mention only that Voevodsky hierarchy does not
{\it replace} sets but rather systematically embeds set--theoretical and categorical constructions and intuitions
into a vaster universe where continuous and discrete are treated on an equal footing.

\newpage

\centerline{\bf 3. Language--like mathematical structures and metamathematics}

\medskip
 {\bf 3.1. Metamathematics: mathematical studies of formalized languages of mathematics.}
Philosophy of mathematics   in the XX--th century had to deal with lessons of
{\it metamathematics,}  especially of G\"odel's incompleteness theorem.
\smallskip

As I have already said, I will consider metamathematics as {\it a special chapter of mathematics itself},
whose subject is the study of {\it formal languages and their interpretations.} On the foreground here 
were the first order formal languages, a formalization
of Euclid's and Aristotle's legacy.  Roughly speaking, to Euclid we owe the mathematics
of spatial imagination (and/or kinematics of solid bodies), whereas  Aristotle founded
{\it the mathematics of logical deduction}, expressed in ``Elements''  by natural language
{\it and} creative usage of drawings.

\smallskip

An important parallel development 
of formal languages involved languages formalizing {\it programs for and processes of computation,}
of which chronologically first in the XX--th century was {\it Church's lambda calculus. }
\smallskip

An important feature of lambda calculus  is the absence of formal distinctions between
the language of programs and the language of  input/output data (unlike
Turing's machines, where a machine ``is'' the program, whereas input/output
are represented by binary words). When, due to von Neumann's insights, this feature
became implemented in hardware, lambda calculus was rediscovered and
 became in the 1960's the basis of development
of programming languages.
\smallskip

These languages are {\it linear}, in the
following sense: the set of all syntactically correct expressions in a formal language
$L$ could be described as a Bourbaki structure consisting of a certain {\it words,} -- finite sequences of letters
 in a given {\it alphabet}, and finite sequences of such words, {\it expressions}.
Words and expressions must be {\it syntactically correct} (precise description of this is a part of
definition of each concrete language). Letters of alphabet 
are subdivided into {\it types}: variables, connectives and quantifiers, symbols for operations, relations ...
Syntactically correct expressions can be terms, formulas, ...

\smallskip
Such Bourbaki structures can be sufficiently rich to produce formal versions of real mathematical texts,
 existing and potential ones,
and make them an object of study.

 \smallskip
 
 I will explain how the advent of category theory (and, to a lesser degree, theory of computability)
 required enriched languages, that after formalization become at first {\it non--linear,}
 and then {\it multidimensional.} Such languages  require for their study  {\it homotopy theory}
 and suggest a respective enrichment of the universe in which interpretations/models are
 supposed to live,
 from  Sets to Homotopy Types, as in the Voevodsky's project (cf. above).

\medskip

 {\bf 3.2. One--dimensional languages of diagrams and graphs.} 
With the development of homological algebra and category theory in
the second half of the XX--th century, the language of commutative
diagrams began to penetrate ever wider realms of mathematics. It took some time
for mathematicians to get used to ``diagram-chasing.'' A simply looking algebraic identity
$kg=hf$, when it expresses a property of four morphisms in a category,
means that we are conemplating a simple {\it  commutative diagram},
in which, besides morphisms $f,g,k,h$, also the objects $A,B,C,D$ invisible in the 
formula $kg=hf$ play key roles:

$$
\xymatrix{A \ar[d]_{f} \ar[r]^{g} &  B
\ar[d]^k \\
C \ar[r]^h & D}
$$

Although this square is not a ``linear expression'', one may argue that it, and its
various generalizations of growing size (even the whole relevant category),
 are still ``one--dimensional''. This means that 
they can be encoded in a graph, whose  vertices are labeled by (names of)
objects of our category, whereas edges are labeled by pairs consisting
of an orientation and a morphism between the relevant objects.

\smallskip

Similarly, a program written in a linear programming language can be encoded in a graph whose 
vertices are labeled by (names of) elementary operations that can be performed
over the relevant data. To understand labeling of (oriented) edges, one must imagine that they
encode channels, forwarding output data calculated by the operation at the start (input) of the edge to the
its endpoint where they become input of the next operation (or the final output, if the relevant vertex is labeled respectively). Labels of edges might then include {\it types} of the relevant data.

\medskip
 {\bf 3.3. From graphs to higher dimensions.} Generally,  a square of morphisms as above 
need not be commutative
(i.~e.~ it is possible that $kg\ne hf$). In order to distinguish these two cases graphically,
we may decide to associate with {\it a commutative square  the  two--dimensional picture}, by glueing 
the interior part of the  square to the relevant graph. 

\smallskip

A well known generalization of this class of spaces are {\it cell complexes}, or, in more
combinatorial and therefore more language--like version, {\it simplicial complexes.}
Of course, we must allow  labels of cells as additional structures.

\smallskip

In this way, we can get, for example, a geometric encoding 
of the category $\Cal{C}$ by a simplicial complex,
in which labeled $n$--complexes are sequences of morphisms
$$
\xymatrix{X_0 \ar[r]^{f_0}  & X_1 \ar[r]^{f_1} & ...   \ar[r]^{f_{n-1}} & X_n}
$$
whereas the face map $\partial^i$ omits one of the objects $X_i$ and, if $1\le i\le n-1$,
replaces the  pair of arrows around $X_i$ by one arrow labeled by the composition of the relevant
morphisms. The resulting simplicial space encodes the whole category
in a simplicial complex that is called {\it the nerve} of the category. Clearly,
not only objects and morphisms, but also all compositions of morphisms 
and relations between them can be read off it.

\smallskip

Thus the language of commutative diagrams becomes a chapter of algebraic topology,
and when the study of functors is required, the chapter of
homotopical topology.

\medskip

{\bf 3.4. Quillen's homotopical algebra and univalent foundations project.} In his influential book [Qu] Quillen
developed the idea that the natural 
language for homotopy theory should appeal {\it not} to the
initial intuition of continuous deformation itself, but rather to a codified
list of properties of category of topological spaces stressing those that are
relevant for studying homotopy.

\smallskip

Quillen defined {\it a closed model category} as a category endowed
with three special classes of morphisms: {\it fibrations,
cofibrations},  and {\it weak equivalences.} The list of axioms to which these 
three classes of morphisms must satisfy is not long but structurally quite sophisticated.
They can be easily defined in the category of topological spaces using homotopy
intuition but remarkably admit translation into many other situations.
An interesting new preprint [GaHa] even suggests the definition of these
classes in appropriate categories of discrete sets, contributing new
insights to old Cantorian problems of the scale of infinities.

\smallskip

Closed model categories become in particular a language of preference for  many contexts 
in which objects of study are quotients of ``large'' objects by ``large'' equivalence relations,
such as homotopy.

\smallskip

It is thus only natural that the most recent Foundation/Superstructure,
 Voevodsky's Univalent Foundations Project (cf. [Vo] and [Aw]) is based on
 direct axiomatization of the world of homotopy types.
 
 \smallskip
 
As a final touch of modernism, the metalanguage of this project is a version of typed
lambda calculus, because its goal is to develop a tool for the  computer assisted verification
of programs and proofs. Thus computers become more and more involved in the interpersonal habitat
of ``theoretical'' mathematics. 

\smallskip

It remains to hope that humans will not be finally excluded from this habitat, as some aggressive
proponents of databases replacing science suggest (cf. [An]).

\vskip1cm

\centerline{\bf Post Scriptum: Truth and Proof in Mathematics}

\medskip

As I have written in [Ma2] , the notion of ``truth''  in most philosophical contexts is a reification
of a certain relationship between humans and {\it texts/utterances/statements}, the relationship
that is called ``belief'', ``conviction'' or ``faith''. 
\smallskip

Professor Blackburn in [Bl]  in his keynote speech to the Balzan Symposium on ``Truth''
(where [Ma2] was delivered) extensively discussed other relationships of humans
to texts, such as {\it scepticism, conservatism, relativism, deflationism.}
However, in the long range all of them are secondary in the practice of a researcher
in mathematics.

\smallskip
I will only sketch here what must be said about  texts, sources
of conviction, and methods of conviction peculiar to mathematics.

\smallskip

{\it Texts.} Alfred North Whitehead  said that all of 
Western philosophy was but a footnote to Plato.

\smallskip

The underlying metaphor of such a statement is: ``Philosophy is a text'',
the sum total of all philosophic utterances.

\smallskip

Mathematics decidedly is {\it not} a text, at least not in the same
sense as philosophy. There are no authoritative
books or articles  to which subsequent 
generations turn again and again
for wisdom. Already in the XX--th century, researchers did not read Euclid,
Newton, Leibniz or Hilbert in order to study
geometry, calculus or mathematical logic. 
The life span of any contemporary mathematical paper or book
can be years, in the best (and exceptional) case decades.
Mathematical wisdom, if not forgotten, lives as an invariant of all
its (re)presentations in a permanently self--renewing discourse. 

\smallskip
 
{\it Sources and methods  of conviction.} Mathematical truth is not revealed, 
and its acceptance is not imposed by any authority.

\smallskip

Ideally, the truth of a mathematical statement is ensured by {\it a proof},
and the ideal picture of a proof is a sequence of elementary arguments
whose rules of formation are explicitly laid down before the proof even begins,
and ideally are common for all proofs that have been devised
and can be devised in future. The admissible starting points
of proofs, ``axioms'', and terms in which they are formulated,
should also be discussed and made explicit.

\smallskip

This ideal picture is so rigid that it became
the subject of mathematical study in metamathematics.

\smallskip

But in the creative mathematics, the role of proof is in no way restricted to
its function of carrier of conviction. Otherwise, there would be no need
for Carl Friedrich Gauss to consider eight (!) different proofs of 
the quadratic reciprocity law (cf. [QuRL] for an extended bibliography; I am grateful to Prof. Yuri Tschinkel for this
reference).
\smallskip

One metaphor of proof is a route, which might be
a desert track  boring and
unimpressive  until one finally reaches the oasis of one's destination, or
 a foot path in green hills, exciting and energizing, opening great vistas
of unexplored lands and seductive offshoots, leading far away even after 
the initial destination point has been reached.

\bigskip

\bigskip

\centerline{\bf References}

\medskip

[An] Ch.~Anderson. {\it The End of Theory.} in: Wired, 17.06, 2008.
\smallskip
[At] M.~Atiyah. {\it Geometry and Physics of the 20th Century.} In [KFNS], p.~4--9.

\smallskip

[Aw] St.~Awodey. {\it Type theory and homotopy.}  (online publication).

\smallskip
[By] Oliver Byrne. {\it The first six Books of the Elements of Euclid.} Facsimile
of the first edition,  published by Taschen.

\smallskip

[De]  S.~Dehaene. {\it The Number Sense. How the Mind creates Mathematics.}
Oxford UP, 1997.

\smallskip

[GaHa] M.~Gavrilovich, A.~Hasson. {\it Exercices de Style: a homotopy thory for set theory.}
arXiv:1102.5562

\smallskip

[He] R.~Hersh. {\it Wings, not foundations!}  In: Essays on Foundations of Mathematics and Logic,
Polimetrica Int.~Sci.~Publisher, Milan, Italy, 2005.
\smallskip

[Ju] A.~Jung. {\it A short introduction to the
Lambda Calculus.} (online publication).
\smallskip
[KFNS] {\it G\'eom\'etrie au XXe si\`ecle. Histoire et Horizons.} (Ed. J.~Kouneiher, D.~Flament,
Ph.~Nabonnand, J.-J.~Szczeciniarz.)
Hermann, Paris, 2005.
\smallskip 
[Ma1] Yu.~Manin. {\it A Course of Mathematical Logic for Mathematicians.}
2nd Edition, with new Chapters written by Yu.~Manin and B.~Zilber. Springer, 2010.

\smallskip

[Ma2] Yu.~Manin.  {\it Truth as value and duty: lessons of mathematics. }
In: Truth in Science, the Humanities, and Religion. Balzan Symposium 2008.
Ed. by N.~Mout and W.~Stauffacher. Springer, 2010, pp. 37--45.
Preprint math.GM/0805.4057. 

\smallskip
[Qu] D.~Quillen. {\it Homotopical Algebra.} Springer LNM, vol. 43, Berlin, 1967.

\smallskip

[QuRL]  {\it Proofs of the Quadratic  Reciprocity Law.}

 http://www.rzuser.uni-heidelberg.de/~hb3/rchrono.html

\smallskip

[Vo] V.~Voevodsky. {\it Univalent foundations project.} Oct.~1, 2010. (online publication).

\bigskip

\enddocument